# Energy Management System for a Low Voltage Direct Current Microgrid: Modeling and experimental validation

Yanandlall Gopee
*LAAS-CNRS, Université de Toulouse, Université Toulouse III - Paul Sabatier*
Toulouse, France
ygopee@laas.fr

Margot Gaetani-Liseo
*LAAS-CNRS, Université de Toulouse, CNRS*
Toulouse, France
mgaetani@laas.fr

Anne Blavette
*SATIE, Univ. Rennes, ENS Rennes CNRS*
Rennes, France
anne.blavette@ens-rennes.fr

Guy Camilleri
*IRIT, Université de Toulouse, Université Toulouse III - Paul Sabatier*
Toulouse, France
Guy.Camilleri@irit.fr

Xavier Roboam
*LAPLACE, Université de Toulouse CNRS*
Toulouse, France
xavier.roboam@laplace.univ-tlse.fr

Corinne Alonso
*LAAS-CNRS, Université de Toulouse, Université Toulouse III - Paul Sabatier*
Toulouse, France
alonsoc@laas.fr

***Abstract*—In the field of microgrids with a significant integration of Renewable Energy Sources, the efficient and practical power storage systems requirement is causing DC microgrids to gain increasing attention. However, uncertainties in power generation and load consumption along with the fluctuations of electricity prices require the design of a reliable control architecture and a robust energy management system for enhancing the power quality and its sustainability, while minimizing the associated costs. This paper presents a mixed approach illustrating both simulation and experimental results of a grid-connected DC microgrid which includes a photovoltaic power source and a battery storage system. Special emphasis is placed on the minimization of the total operating cost of the microgrid while considering the battery degradation cost and the electricity tariff. Thereby, an optimal energy management system is proposed for Energy Storage Systems scheduling and enabling the minimization of the electricity bill based on simple models. Simultaneously, the differences between simulation and laboratory performances are highlighted.**

***Keywords—Microgrid, energy management system, optimization, linear programming, battery degradation***

## I. Introduction

With a considerable increase of Renewable Energy Sources (RESs) in electrical energy systems, the concept of Microgrids (MGs) has emerged a novel paradigm for active distribution networks contributing to the emergence of smart grids. The latter concept has been gaining increasing interest by the research community during the last 20 years [1]. MGs can be connected to the main electricity distribution grid in order to improve reliability of power production and reduce power losses due to the proximity between energy productions and consumptions as shown in Fig. 1. Nevertheless, they can also be operated in an islanded mode to improve the resilience, ensuring the power supply continuity in case of grid failure as they are potential solutions for electrification of off-grid islands, rural villages, and remote areas [2], [3]. Furthermore, in case of default of the utility grid, some MGs can collaborate with other MGs in the neighborhood. Nowadays, MGs can be fine-tuned so as to minimize their associated energy bills, while consuming more renewable energy and be better integrated in their environment. Moreover, technical adjustments enable to mitigate the control and coordination problems which can originate from inter-temporal operation of Energy Storage Systems (ESSs) and intermittent nature of RESs issues. Thus, an Energy Management System (EMS) becomes indispensable in order to explicitly handle constraints, energetics and economic criteria thereby enabling optimal scheduling of power resources and ESSs in MGs so as to achieve optimal supply-demand balance [4]. EMS is the tertiary layer of MG control operating hierarchically above the real-time control and the power management levels, which are the primary and secondary layers respectively [5].

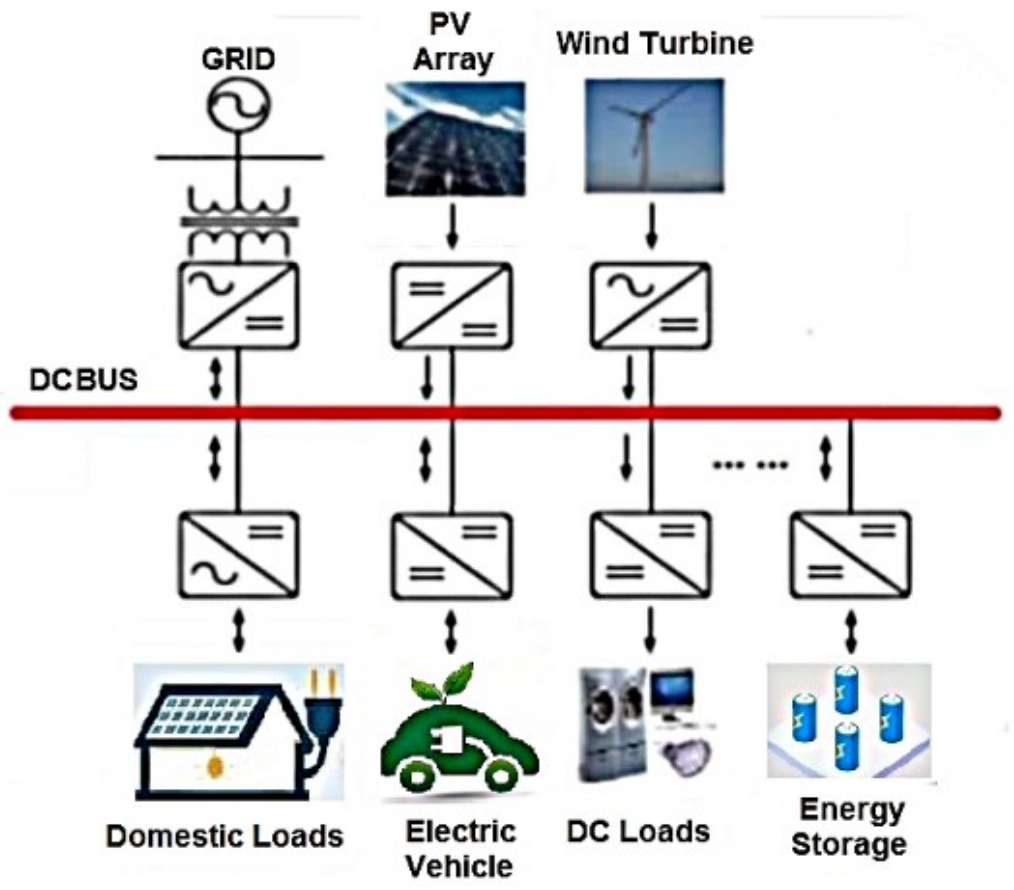

Fig. 1. Schematic diagram of a generic MG.

In the recent years, a number of research papers related to EMS dedicated to MGs have been published. In [6], a review of several optimization methods of energy management is presented. A controller for EMS called the Resilient Model Predictive Controller (RMPC) based on the Model Predictive Control (MPC) principles is developed in [7] for grid connected residential MG in order to minimize the electricity bill and load shedding. Furthermore, in [8] a Resilient Open Loop Feedback Controller (ROLFC) for EMS is developed using a multistage stochastic optimization (MSO) approach. In [9], a Linear Programming (LP) along with heuristic algorithm is proposed for minimization of operating cost of a grid connected MG. In [10], a heuristic approach for energy management of an alternating current (AC) MG is proposed which results in non-optimal solutions. A multi-objective predictive control strategy is proposed in [11] for grid-connected photovoltaic (PV) and battery systems for residential buildings aiming to dynamically decrease the cost of electricity compared to traditional systems while reducing the $CO_2$ emissions. In [12], an economic cost model is used for optimal management of MG Distributed Energy Resources (DERs), taking into account the degradation cost of the cyclic ageing of the battery. However, these EMSs have not been experimentally validated on facilities with actual devices or based on hardware-in-the-loop platforms.

Concerning experimentally verified EMSs, an EMS for battery-based hybrid grid-connected MG is designed and experimentally tested in [13] which enables the scheduling of optimal power references for DERs and also allows to reduce cost in the MG over a 24-hour ahead forecast data. In [14], a Mixed-Integer Linear Programming (MILP) mathematical model is proposed in order to optimally schedule the power references of the DERs in a grid-connected hybrid PV-wind-battery MG. In [15], an EMS is developed that optimally balances generation and demand with scheduled power generation of PV, wind, and tidal energy sources with the objectives of maximizing the utilization of RES for an islanded Direct Current (DC) MG. With reference to [15], the battery ageing cost is considered in detail taking into account several parameters such as temperature, operation and maintenance cost and residual value which the author obtained from the practical regression models of battery degradation developed in [16].

Moreover, a dynamic optimization model is proposed in [17] for power sharing among sources in an islanded DC MGs, but the battery degradation cost for sustainable battery lifetime is not taken into consideration. In [18], a MPC approach is applied for achieving economic efficiency in a grid-connected MG operation management. Thereby, three different experiments are executed with different planning horizons in order to assess the feasibility of the MPC–MILP control scheme. A predictive EMS is developed in [19] for an islanded PV-battery DC MG for preventing or reducing outage duration. In [20], a distributed control and non-optimal EMS is proposed for preventing State of Charge (SoC) violations of hybrid ESSs in islanded DC MGs, however optimization model for cost-effective MG operation is not included.

Most of the aforementioned EMSs have focused on economic optimization and have been experimentally validated using hardware-in-the-loop platforms. However, the large majority of researches in the existing literature have ignored the battery degradation cost in the latter models. Nevertheless, relatively few studies [9], [12], [15], [16], [18] have considered the wear and tear cost of the battery. Moreover, the latter previous works have not explored the difference between simulation and experimental results. As a substantive contribution in the field of power distribution in MGs, this research not only considers the EMS for maximizing economic benefit and optimally scheduling of ESS but it also considers the battery degradation cost and investigates the discrepancy between prototype and real-life performance. Furthermore, this paper explicitly includes some details on the implementation architecture.

In summary, considering all of the aforementioned limitations regarding previous studies on EMS of MG, the following contributions are made in this research:

- An EMS strategy is developed, which performs economic dispatch for minimizing electricity cost of a grid-connected Low Voltage Direct Current (LVDC) MG while satisfying the user's demand in terms of uninterrupted power supply. It takes into consideration the degradation cost of the ageing cycle of the battery in order to reflect the true operational cost of a MG.
- The proposed EMS strategy is validated by an experimental set-up. Analyses are carried out to highlight the differences between simulation and experimental performances.

This paper is organized into five sections. Section II illustrates the case study of an LVDC grid-connected MG operation and also details the EMS mathematical modeling. Section III describes the set-up of the MG experimental platform. Section IV presents the simulation and experimental results of the proposed EMS strategy and finally, Section V presents the conclusion of the paper.

## II. CASE STUDY

### *A. System description*

A grid-connected residential PV-battery LVDC MG with a dc-coupled architecture is studied in order to enable experimental validation, as illustrated in Fig. 2. It consists of a discrete time model of a PV system which includes a storage capacity for self-consumption by a user who is connected to the conventional grid. The system environment includes two input variables $P_{pv}$ which is the solar production potential and $P_{load}$ which is the power demand. The decision variable ($P_{st}$) is the power for charging/discharging the battery, the recourse variable ($P_g$) is the power exchange with the grid and ($P_{batt}$) is the internal power of the battery.

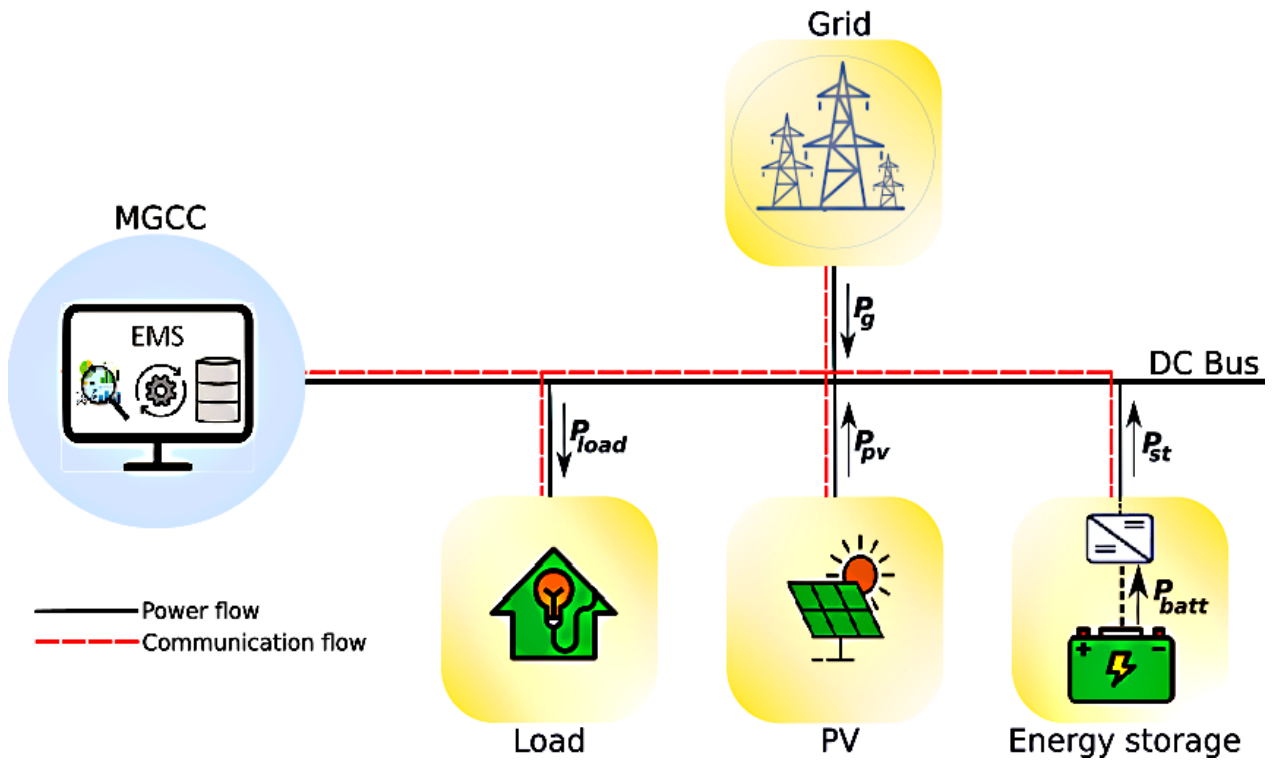

Fig. 2. Grid connected DC MG architecture.

An EMS is used for supervisory control and a day-ahead optimal power scheduling of the LVDC MG. The Microgrid Central Controller (MGCC) computes the inflow and outflow of the power profile offline, for each flexible entity, namely the storage units and the grid, in order to minimize the cost of electricity over a given day and to ensure power balance in the MG. It is assumed that the model is a deterministic one where $P_{pv}$ and $P_{load}$ of the household are already known in advance and also that the conventional grid is always available.

### *B. Optimization problem*

The optimization problem for achieving the techno-economical goals of the grid connected LVDC MG is mathematically modeled using LP formulation as in [21] instead of MILP algorithms. This choice of a LP in the control allows to face both storage and grid management constraints such as preventing simultaneous charging/discharging of the battery and purchase/selling of energy to/from the grid respectively which are guaranteed cost-wise without using binary variables. The latter variables would have transformed the problem into a MILP formulation that would involve using MILP solvers instead of LP solvers which may result in a longer computing time for larger-scale examples. Hence, in this study LP was preferred in the perspective of future work.

All the variables at the bus level of this generic minimization problem are linked through the power balance relationship at time t which is given as follows:

$$P_g(t) + P_{st}(t) + P_{pv}(t) = P_{load}(t) \quad (1)$$

The recourse variable $P_g$ and the decision variable $P_{st}$ are separated into two components to fulfill the linearization of the system:

$$P_g(t) = P_g^+(t) + P_g^-(t) \quad (2)$$

$$P_{st}(t) = P_{st}^{dis}(t) + P_{st}^{ch}(t) \quad (3)$$

$P_g^+(t)$ is the power drawn from the grid whereas $P_g^-(t)$ is the power sold to the grid. $P_{st}^{dis}(t)$ and $P_{st}^{ch}(t)$ are the discharging and charging power of the battery respectively.

*1) Objective function and Constraints*

The objective function that minimizes the cost of the electricity exchanged with the grid is given by (4). $C_E$ is the electricity cost to be minimized over the time horizon $t_{day}$, which is sampled at the time interval $\Delta t$.

$$min\ C_E = min\ \Delta t \sum_t^{t_{day}} \left[c_g^+(t)P_g^+(t) + c_g^- P_g^-(t) + c_{st}P_{st}^{dis}(t) + c_{st}P_{st}^{ch}(t)\right] \quad (4)$$

Where $c_g^+(t)$ is the electricity purchase tariff at time t which can vary according to off- or on-peak energy prices. The costs which are constant over the time horizon are: $c_g^-$, the selling price of energy to the grid and $c_{st}$, the battery ageing unit cost associated with the charging/discharging decision variables.

In parallel to (1), the following constraints must be satisfied at time t in order to achieve a feasible solution:

$$P_{gmax}^- \leq P_g(t) \leq P_{gmax}^+ \quad (5)$$

Where $P_{gmax}^+ \geq 0$ and $P_{gmax}^- \leq 0$ are defined by the contracted powers subscribed by the user.

Consequently, $P_{batt}(t)$ is calculated linearly as follows:

$$P_{batt}(t) = \frac{P_{st}^{dis}(t)}{\eta_{cvs}} + P_{st}^{ch}(t).\eta_{cvs}.\eta_E \quad (6)$$

Where $\eta_{cvs}$ is the efficiency of the converter and $\eta_E$ is the energy efficiency of the battery which is more precisely its round-trip efficiency considered at charging mode only. Thus, the boundary of the charging/discharging powers is:

$$P_{battchargemax} \leq P_{batt}(t) \leq P_{battdischargemax} \quad (7)$$

Where $P_{battdischargemax} \geq 0$ and $P_{battchargemax} \leq 0$ are related to the battery C-rate which is the measurement of current in which a battery is charged and discharged. In addition, $P_{st}^{dis}(t)$ and $P_{st}^{ch}(t)$ are also bounded:

$$0 \leq P_{st}^{dis}(t) \leq P_{battdischargemax}.\eta_{cvs} \quad (8)$$

$$\frac{P_{battchargemax}}{\eta_{cvs}.\eta_E} \leq P_{st}^{ch}(t) \leq 0 \quad (9)$$

The power exchange with the battery is indirectly constrained through the energy stored in the ESS allowing the calculation of the State of Energy (SoE). The remaining and available energy of a battery is defined by the SoE indicator with reference to its nominal energy capacity and compared to the SoC which estimates what is left in the battery in terms of the amount of charge related to its nominal capacity [22]. The SoE is calculated as given below:

$$SoE(t+1) = SoE(t) - \frac{P_{st}^{dis}(t)\times\Delta t}{\eta_{cvs}\times E_{nom}} - \frac{P_{st}^{ch}(t).\eta_{cvs}.\eta_E\times\Delta t}{E_{nom}} \quad (10)$$

Where $E_{nom}$ is the nominal energy capacity of the battery. However, the battery has a limited storage capacity:

$$SoE_{min} \leq SoE(t) \leq SoE_{max} \quad (11)$$

*2) Ageing cost of the battery and grid cost*

The battery ageing unit cost ($c_{st}$) is expressed as the cost of exchanging 1 kWh. Therefore, the total amount of energy that can be exchanged by the battery over its lifetime ($E_{lifebatt}$), as derived from [23], is calculated as follows:

$$E_{lifebatt} = E_{nom} \times N_{cycles} \times DoD \quad (12)$$

Where $N_{cycles}$ is the number of cycles and $DoD$ the Depth of Discharge of the battery. Knowing the battery's price ($c_{batt}$), $c_{st}$ can be calculated by the following equation:

$$c_{st} = \frac{1}{2}\frac{c_{batt}}{E_{lifebatt}} \quad (13)$$

A factor of $\frac{1}{2}$ is used to account for both charging and discharging modes during a cycle [16]. This is a simple ageing model and the battery price is estimated at *85 €/kWh*, leading to a cost of $c_{st}$ equal to *0.235 €/kWh*. Moreover, according to the power subscription by the user and based on off- or on-peak energy price in France (regulated tariffs), the actual purchase price of energy $c_g^+(t)$ from the grid is *0.1360 €/kWh* during off-peak hours and *0.1821 €/kWh* during on-peak ones. The selling price ($c_g^-$) is constant and equal to *0.10 €/kWh* for a PV system ≤ 9 kWc. Nevertheless, as the electricity price has been increasing in the recent years, based on the evolution curve of the electricity tariff [24], [25], the hypothesis made in this study is that in the future, the purchase tariff for electrical power will be five-fold higher than the actual price and the selling price will increase only two-fold. This extrapolation is clearly hypothetical but the recent cost evolution of energy prices, especially the EPEX-SPOT tariff with selling prices of electricity varying from *0* to *1000 €/MWh* makes credible this assumption [25]. Therefore, the purchase electricity costs are considered as *0.68 €/kWh* during off-peak and *0.9105 €/kWh* during on-peak hours. The selling electricity rate is considered as *0.20 €/kWh*. It is to be noted that off-peak hours can be consecutive or scattered within several time slots. There are 8 off-peak hours per day which are concentrated at different times of the day when demand is naturally less sustained. The 8 off-peak hours per day considered in this study is comprised between 12 am and 8 am.

## III. Set-up Description of the Experimental Microgrid

A test bench is developed at the "Laboratoire d'Analyse et d'Architecture des Systèmes" (LAAS) to emulate the operation of a laboratory scale LVDC MG platform as shown in Fig. 3. The latter bench includes a programmable power supply and a loading facility where both consist of integrated converters, thus enabling the emulation of PV production and load consumption respectively. The test bench consists of an ESS which is connected to the DC bus via a bidirectional current converter. Thus, the ESS and the associated converter can operate under real-life conditions involving actual constraints such as variation in the efficiency of the ESS and the converter, change of temperature at the location, etc.

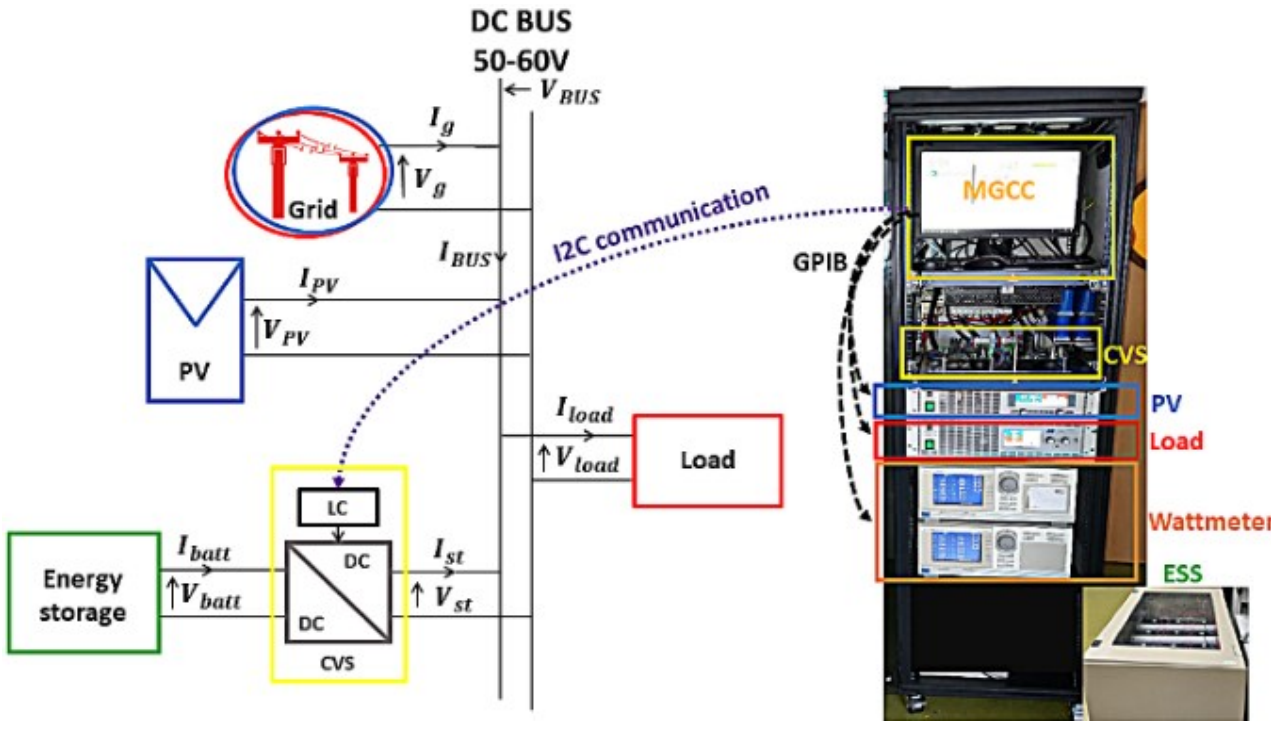

Fig. 3. Experiment platform of the LVDC MG.

The MG is controlled using the DC Bus Signaling (DBS) which is a control strategy where each entity controls the bus voltage level according to the energy exchange in the LVDC MG in order to maintain its power balance [26]. Additionally, the above-mentioned bench is driven following a hierarchical control. The primary control is a decentralized structure used to control the converter as has been developed at LAAS [27]. The Local Controller (LC) present on the converter board ensures the above-mentioned primary control. The secondary and tertiary controls are centralized structures incorporated for ensuring the functions of Power Management System (PMS) by the DBS and the EMS respectively. The MGCC communicates via an Inter Integrated Circuit ($I^2C$) communication bus with the LC for the converter associated with the ESS. The PV, load and bench instrumentation system communicates through General Purpose Interface Bus (GPIB) communication protocol.

In this study, a pack of lead–acid battery is connected as the ESS of the MG. The battery rated energy capacity is 2.16 kWh and the maximum power for charging and discharging is 0.54 kW at 0.25 C-rate and 0.846 kW at 0.4 C-rate respectively [28]. The energy efficiency of the battery ($\eta_E$) is 83% and the converter efficiency ($\eta_{cvs}$) is 96% [27]. The maximum and minimum limits of the battery SoE are defined as 95% and 5% respectively. The maximum power of $P_g$ that can be drawn or transferred to the grid is taken as 0.9 kW. The voltage, current and power limit values of the DC bus are given below in Table I.

TABLE I. GENERAL CHARACTERISTICS OF THE LVDC MG

| Parameter | Values and Units |
|---|---|
| [$V^{bus}_{min}$ ; $V^{bus}_{max}$] | [50 ; 60] V |
| [$I^{bus}_{min}$ ; $I^{bus}_{max}$] | [-20 ; 20] A |
| $P^{bus}_{max}$ | 1 kW |

## IV. RESULTS AND DISCUSSION

The PV production and load consumption data for a 24-hour ahead period used for the simulation and experimental set-up is extracted from the "Oahu Archive web site" and "IEEE PES Test Feeder" respectively. These latter values are scaled for the purpose of experimental validation. Fig. 4 shows the PV production, load consumption and net power ($P_{net}$) which is the power difference between $P_{pv}$ and $P_{load}$ as expressed hereunder:

$$P_{net} = P_{pv} - P_{load} \tag{14}$$

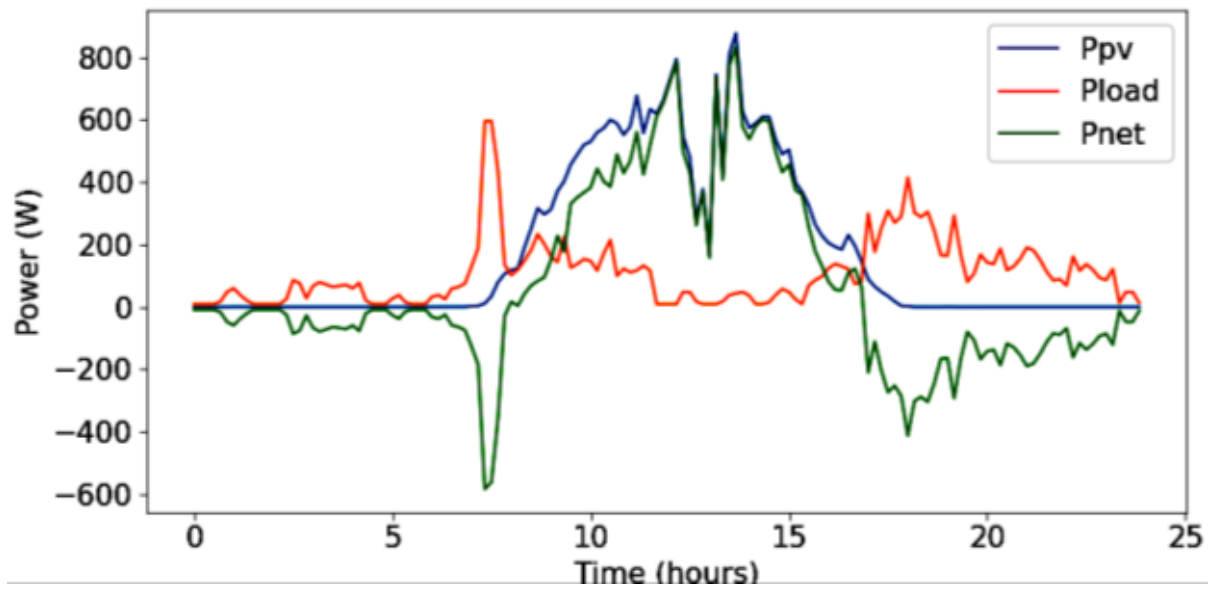

Fig. 4. PV production, load consumption and net power.

### A. Simulation results

The simulation is carried out over a 24-hour prediction horizon. The power references for a LP based optimal scheduling are computed offline for a sampling rate of 10 minutes (Δt = 10 minutes) by the MGCC. The model is developed in Python using the pyomo modeling language [29] and solved using the GUROBI solver. The proposed EMS is tested with an initial battery SoE maintained at 35% as shown in Fig. 5. It is underlined that for the simulation performed the global optimal solution is achieved, as described below.

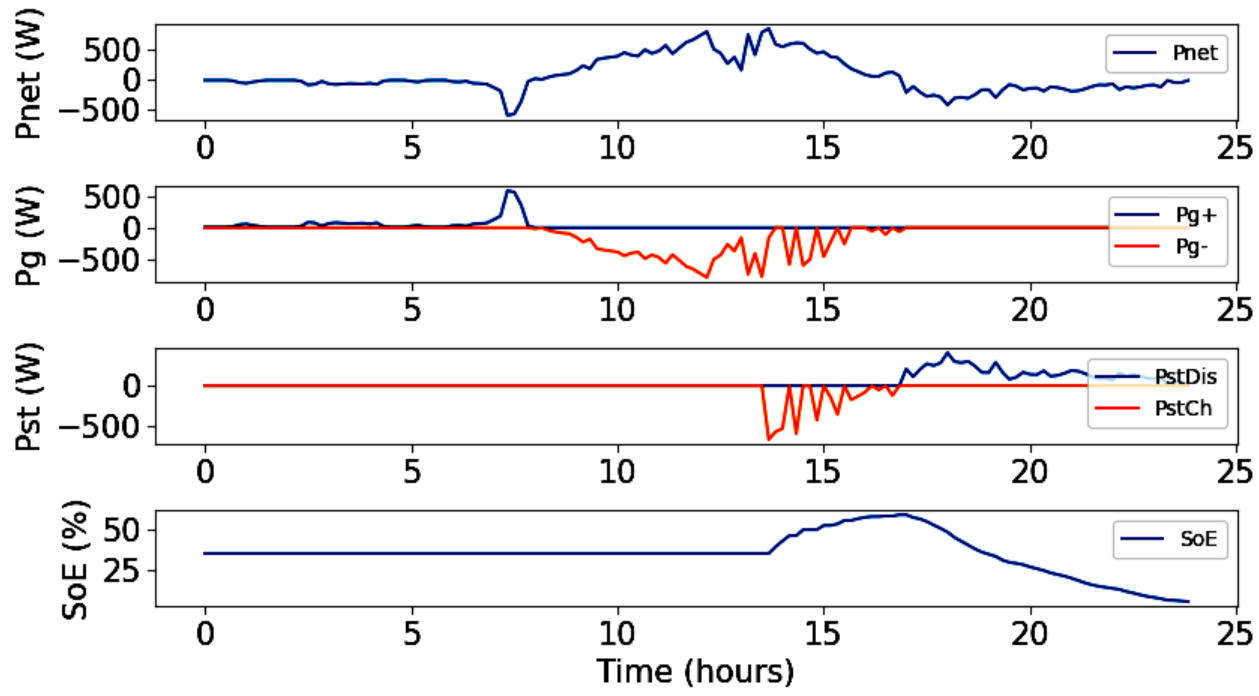

Fig. 5. Simulation results for an initial SoE of 35%.

The optimal solution for a minimal cost is achieved by enabling the utilization of the grid during off-peak hours between 12 am and 8 am. The battery is charged between 1:30 pm and 4:50 pm when $P_{net} > 0$ starting with an initial SoE of 35% to reach a required level of SoE of about 59% so as to be able to meet the user's energy demand during peak hours without drawing power from the grid. It is pointed out that the battery is not charged to its maximum SoE as it is costly for charging and discharging. Also, an SoE of 59% is enough for the time horizon considered. Furthermore, the battery is discharged at the end of the day between 4:50 pm and 12 am so as to avoid drawing power from the grid during low PV production and peak hours. Consequently, the electricity bill is minimized while simultaneously considering the operating cost due to the battery degradation. The electricity cost over the time horizon for this simulation scenario is 0.27€ as per the objective function value. Thus, if the latter value is positive, the user need to pay the electricity bill (deficit) whereas, if the value is negative the user is on profit.

### B. Experimental results

The proposed EMS has been experimentally tested at the lab-scale MG of LAAS. Fig. 6 shows the experimental results of the previous scenario for an initial SoE of 35% of the battery. The latter results allow to assess the validity of the offline optimization model used in simulation. Thus, the optimal trajectories obtained offline based on the LP and on the simple models of PV and battery devices are directly injected as the power references in the actual storage system. The PV production and load consumption are also supplied as power references to the test bench for experimental purposes.

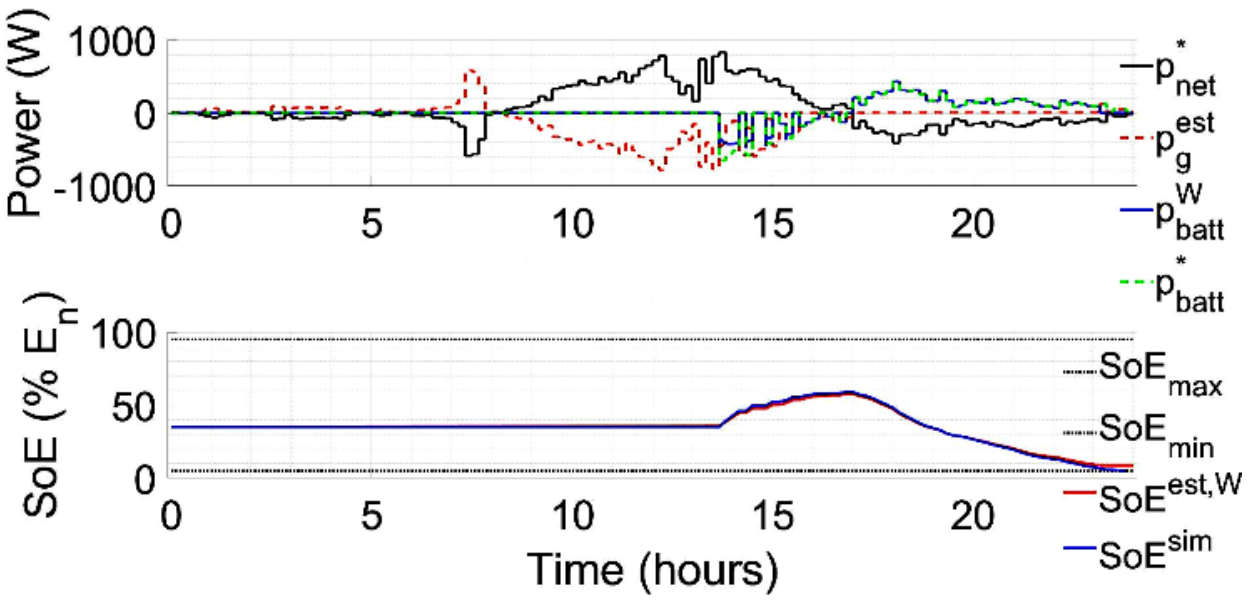

Fig. 6. Experimental results for an initial SoE of 35%.

In Fig. 6, $P^*_{net}$ is the power difference between the PV production and load consumption, $P^{est}_g$ is the estimated power that is drawn and transmitted to the utility grid related to experimental measurement, $SoE_{max}$= 95% and $SoE_{min}$= 5% are the maximum and minimum limits of the SoE of the battery respectively and $SoE^{sim}$ is the resultant SoE of the battery with respect to the offline optimization during simulation. Regarding the experimentation process, $P^W_{batt}$ is the charging/discharging power of the battery measured by the wattmeters enabling the estimation of the battery's $SoE^{est,\ W}$.

Fig. 7 is a close-up of the experimental results between 1 pm and 12 am. It can be noticed that the experimental results slightly differ from the simulation results, especially the charging/discharging power of the battery. It can be specifically observed that between 1:30 pm and 3:30 pm, the power reference $P^*_{batt}$ and the measured power $P^W_{batt}$ are not the same. It is noted that $P^W_{batt}$ is unable to follow the power references of $P^*_{batt}$ during the charging of the battery. This is because the latter had reached its maximum voltage, thereby requiring to be charged at a constant voltage while simultaneously limiting its current according to the classic Constant Current-Constant Voltage (CC-CV) protocol for battery charging phase [30], [31]. Therefore, experimentally the battery is unable to charge as much as expected during simulation. It is also noted that during the experiment process, in the time range between 11 pm and 12 am, the battery is not fully discharged contrarily to simulation results. This is due to the fact that the actual battery is not charged adequately during the experimental process, contrarily to the charge process forecasted during the simulation. Thus, a proportion of the required energy must be drawn from the grid during on-peak hours to meet the user's demand, causing an increase in the electricity bill. Thereby, for the identical scenario as in simulation, the electricity cost over the time horizon for the experimental set-up is 0.33€. Thus, the percentage error of the operation cost between simulation and experimental is 22.2%.

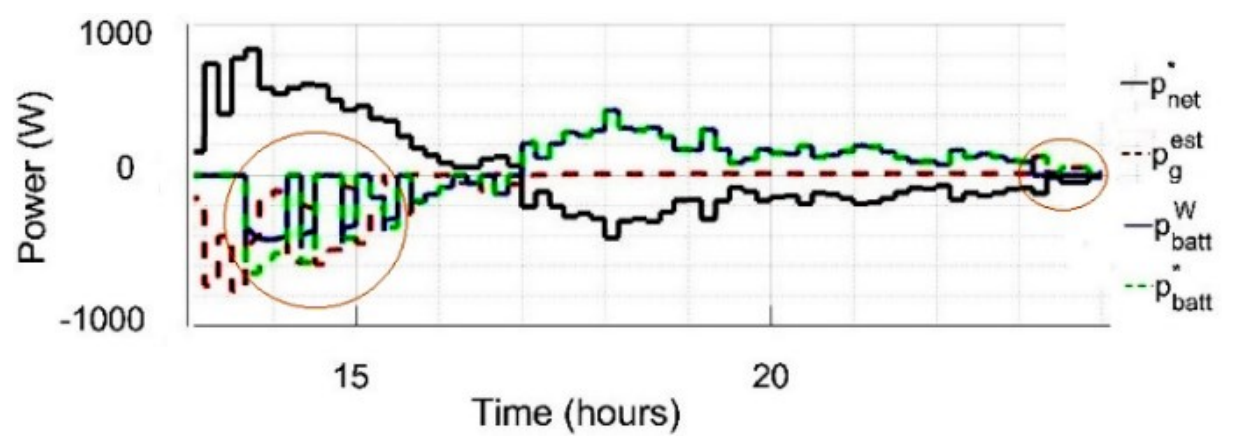

Fig. 7. Close-up of the experimental results for an initial SoE of 35%.

The above-mentioned discrepancy between laboratory experimentation and simulation can be explained by the lack of accuracy of the battery model. In this study, a simple Power-Energy (PE) model is considered for the battery modeling in conformity with most of the EMSs considering an optimization approach as outlined in the literature. This is especially because a linear model derivation must be achieved for LP. Nevertheless, in the existing battery literature, the SoC is the most common indicator used as it gives an accurate estimate of the battery state [31], [32]. Thus, SoC is mainly applied in the field of battery, whereas SoE is mostly employed in EMS optimization approach. However, SoC and SoE are significantly different indicators. Thus, in order to increase accuracy of the battery state, $SoC^{est,W}$ is estimated in the test bench considering the current and the nominal capacity $C_{nom}$ of the battery. On the other hand, the $SoE^{est,\ W}$ is estimated with respect to the current, the voltage and the nominal energy of the battery.

The SoC gives the state of the battery in terms of capacity expressed in ampere hour (Ah) current-wise. Contrarily, the SoE gives the state of the battery in terms of energy expressed in watt-hour (Wh) power-wise. Thus, as neither the current nor the voltage of the battery are considered in the offline optimization model, the SoC cannot be calculated in simulation. Therefore, $SoC^{est,W}$ and $SoE^{est,\ W}$ are estimated [22], [31]–[33] using the values from the wattmeters in the above-mentioned experimental set-up, as follows:

$$SoC^{est,W}(t+1) = SoC^{est,W}(t) - \frac{I^{ch}_{batt}(t).\eta_F \times \Delta t}{C_{nom}} - \frac{I^{dis}_{batt}(t) \times \Delta t}{C_{nom}} \quad (15)$$

$$SoE^{est,\ W}(t+1) = SoE^{est,\ W}(t) - \frac{I^{ch}_{batt}(t) \times V_{batt}(t).\eta_F \times \Delta t}{E_{nom}} - \frac{I^{dis}_{batt}(t) \times V_{batt}(t) \times \Delta t}{E_{nom}} \quad (16)$$

The faradic efficiency of the battery ($\eta_F$) is 96%; $I^{ch}_{batt} \leq 0$ and $I^{dis}_{batt} \geq 0$ are the charging and discharging currents of the battery respectively; $V_{batt}$ is the voltage of the battery.

Fig. 8 shows the difference between the SoE of the battery in simulation and experimental results for the initial scenario of an SoE of 35%. The difference between SoE and SoC can also be observed. $SoC_{max}$ and $SoC_{min}$ are the maximum and minimum limits of the SoC for the battery respectively. $SoC^{est,W}$ is the estimated SoC of the battery during the experimentation process.

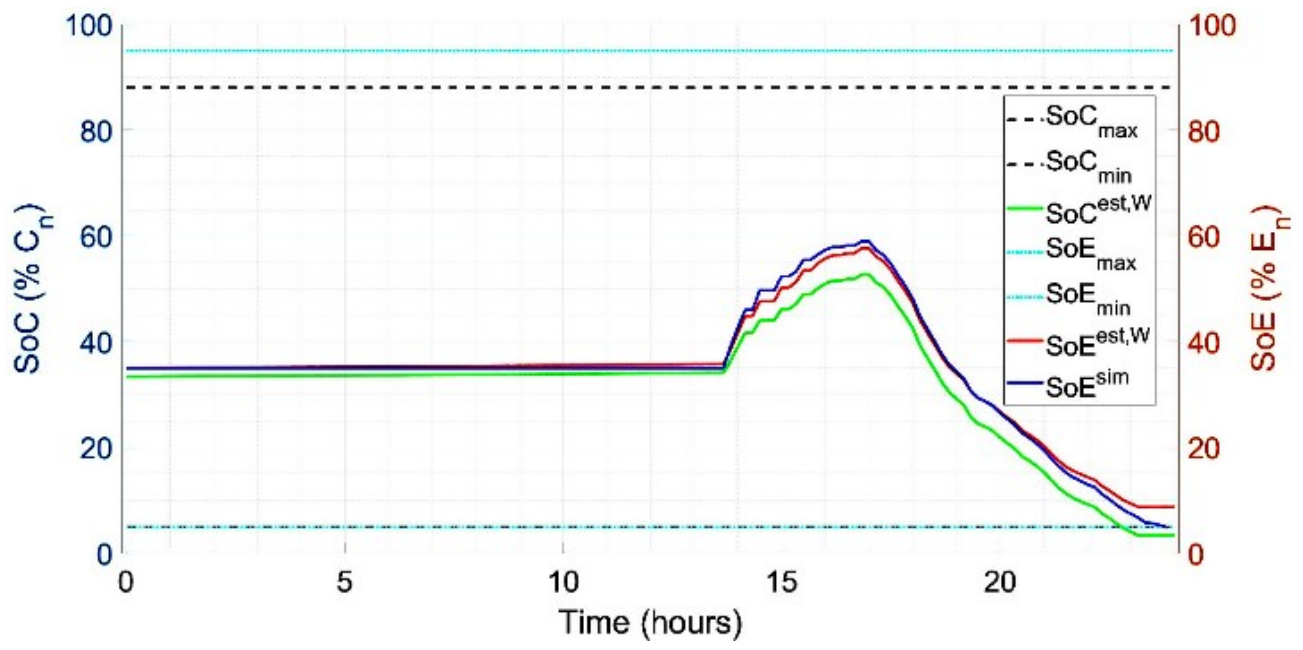

Fig. 8. Comparison of the variations of SoE of 35% between simulation and experimental results relative to SoC limiting values.

As shown in Fig. 8, the battery is unable to be discharged at the end of the process as it has reached its minimum SoC limit. It can also be observed that the approximate deviation of the SoE between simulation and experimental results is around 3.5% at the end of the test. The latter results and the percentage error of the operation cost confirm the major drawback of the PE battery model. Nevertheless, as per the results of this study, it can be concluded that the errors mainly due to the PE model of the battery are not significant over a 24-hour time horizon. However, the accumulation of the aforementioned errors may have a significant impact for optimization studies based on a longer time horizon. In order to minimize this error, a more accurate battery model can be used which consider both the current and the voltage of the battery in the formulation of the optimization model. Furthermore, for the simulation and the optimization process of the MGCC, the SoE value must be updated in real time with reference to the monitoring system or by using a Battery Management System (BMS).

## V. CONCLUSION

In this study, the total operating cost of a grid-connected DC MG is minimized by an EMS strategy. Emphasis is also put on the cost implication with regards to the battery lifetime. Both simulation and experimental studies are carried out in order to validate the performance of the proposed EMS strategy at a laboratory-scale MG. The results confirm the feasibility of the supervisory control and demonstrate the significant impact of the battery model on the cost function and also influence the SoE. It is highlighted that the percentage error between simulation and experimental results is 22.2% for the operation cost and around 3.5% for the deviation of the SoE. Thus, this study confirms the limitation of the PE battery model. This work can be further extended involving a more precise modeling of the battery. In the next stage of this research, the design of a decentralized controller by a Multi-Agent System (MAS) will be carried out for a more complex case study regarding a Multi-MG network. A real-time optimization will also be implemented and further cases with stochastic description of grid failure will be investigated.

## REFERENCES

[1] A. Vasilakis, I. Zafeiratou, D. T. Lagos, and N. D. Hatziargyriou, 'The Evolution of Research in Microgrids Control', *IEEE Open Access J. Power Energy*, vol. 7, pp. 331–343, Oct. 2020.

[2] M. Nasir, H. A. Khan, A. Hussain, L. Mateen, and N. A. Zaffar, 'Solar PV-Based Scalable DC Microgrid for Rural Electrification in Developing Regions', *IEEE Trans. Sustain. Energy*, vol. 9, no. 1, pp. 390–399, Jan. 2018.

[3] P. Nduhuura, M. Garschagen, and A. Zerga, 'Impacts of Electricity Outages in Urban Households in Developing Countries: A Case of Accra, Ghana', *Energies*, vol. 14, no. 12, 3676, Jun. 2021.

[4] D. Wang, J. Qiu, L. Reedman, K. Meng, and L. L. Lai, 'Two-stage energy management for networked microgrids with high renewable penetration', *Appl. Energy*, vol. 226, pp. 39–48, Sep. 2018.

[5] D. Y. Yamashita, I. Vechiu, and J.-P. Gaubert, 'A review of hierarchical control for building microgrids', *Renew. Sustain. Energy Rev.*, vol. 118, 109523, Feb. 2020.

[6] B. Jyoti Saharia, H. Brahma, and N. Sarmah, 'A review of algorithms for control and optimization for energy management of hybrid renewable energy systems', *J. Renew. Sustain. Energy*, vol. 10, no. 5, 053502, Sep. 2018.

[7] J.-J. Prince A., P. Haessig, R. Bourdais, and H. Gueguen, 'Resilience in energy management system: A study case', *IFAC-Pap.*, vol. 52, no. 4, pp. 395–400, 2019.

[8] J.-J. Prince Agbodjan, P. Haessig, R. Bourdais, and H. Gheguen, 'Stochastic modelled grid outage effect on home Energy Management', *IEEE Conference on Control Technology and Applications*, Montreal, QC, Canada, pp. 1080–1085, Aug. 2020.

[9] S. Chakraborty and M. G. Simoes, 'PV-Microgrid Operational Cost Minimization by Neural Forecasting and Heuristic Optimization', in *2008 IEEE Industry Applications Society Annual Meeting*, Edmonton, Alberta, Canada, pp. 1–8, Oct. 2008.

[10] J. B. Almada, R. P. S. Leão, R. F. Sampaio, and G. C. Barroso, 'A centralized and heuristic approach for energy management of an AC microgrid', *Renew. Sustain. Energy Rev.*, pp. 1396–1404, Jul. 2016.

[11] K. Shivam, J.-C. Tzou, and S.-C. Wu, 'A multi-objective predictive energy management strategy for residential grid-connected PV-battery hybrid systems based on machine learning technique', *Energy Convers. Manag.*, vol. 237, 114103, Jun. 2021.

[12] M. A. Hossain, H. R. Pota, S. Squartini, F. Zaman, and K. M. Muttaqi, 'Energy management of community microgrids considering degradation cost of battery', *J. Energy Storage*, vol. 22, pp. 257–269, Apr. 2019.

[13] A. C. Luna, N. L. Diaz, M. Graells, J. C. Vasquez, and J. M. Guerrero, 'Mixed-Integer-Linear-Programming-Based Energy Management System for Hybrid PV-Wind-Battery Microgrids: Modeling, Design, and Experimental Verification', *IEEE Trans. Power Electron.*, vol. 32, no. 4, pp. 2769–2783, Apr. 2017.

[14] A. Luna *et al.*, 'Optimal power scheduling for a grid-connected hybrid PV-wind-battery microgrid system', in *2016 IEEE Applied Power Electronics Conference and Exposition (APEC)*, Long Beach, CA, USA, pp. 1227–1234, Mar. 2016.

[15] M. F. Zia, M. Nasir, E. Elbouchikhi, M. Benbouzid, J. C. Vasquez, and J. M. Guerrero, 'Energy management system for a hybrid PV-Wind-Tidal-Battery-based islanded DC microgrid: Modeling and experimental validation', *Renew. Sustain. Energy Rev.*, vol. 159, 112093, May 2022.

[16] M. F. Zia, E. Elbouchikhi, and M. Benbouzid, 'Optimal operational planning of scalable DC microgrid with demand response, islanding, and battery degradation cost considerations', *Applied Energy*, vol. 237, pp. 695–707, Mar. 2019.

[17] S. Moayedi and A. Davoudi, 'Unifying Distributed Dynamic Optimization and Control of Islanded DC Microgrids', *IEEE Trans. Power Electron.*, vol. 32, no. 3, pp. 2329–2346, Mar. 2017.

[18] A. Parisio, E. Rikos, G. Tzamalis, and L. Glielmo, 'Use of model predictive control for experimental microgrid optimization', *Appl. Energy*, vol. 115, pp. 37–46, Feb. 2014.

[19] D. Michaelson, H. Mahmood, and J. Jiang, 'A Predictive Energy Management System Using Pre-Emptive Load Shedding for Islanded Photovoltaic Microgrids', *IEEE Trans. Ind. Electron.*, vol. 64, no. 7, pp. 5440–5448, Jul. 2017.

[20] J. Xiao, P. Wang, and L. Setyawan, 'Multilevel Energy Management System for Hybridization of Energy Storages in DC Microgrids', *IEEE Trans. Smart Grid*, pp. 847–856, 2015.

[21] R. Sioshansi and A. J. Conejo, *Optimization in Engineering*, vol. 120. Cham: Springer International Publishing, 2017.

[22] K. Mamadou, E. Lemaire, A. Delaille, D. Riu, S. E. Hing, and Y. Bultel, 'Definition of a State-of-Energy Indicator (SoE) for Electrochemical Storage Devices: Application for Energetic Availability Forecasting', *J. Electrochem. Soc.*, vol. 159, no. 8, pp. A1298–A1307, Jul. 2012.

[23] D. Hernández-Torres, C. Turpin, X. Roboam, and B. Sareni, 'Techno-economical optimization of wind power production including lithium and/or hydrogen sizing in the context of the day ahead market in island grids', *Elsevier MATCOM journal*, vol. 158, pp. 162–178, Apr. 2019.

[24] A. Daubaire and J. Pujol, 'Prix de l'énergie : carburants, gaz, électricité, un tiercé en hausse', *Insee*, Dec. 2021.

[25] https://www.rte-france.com/eco2mix/les-donnees-de-marche

[26] J. Schonbergerschonberger, R. Duke, and S. D. Round, 'DC-Bus Signaling: A Distributed Control Strategy for a Hybrid Renewable Nanogrid', *IEEE Trans. Ind. Electron.*, vol. 53, no. 5, pp. 1453–1460, Oct. 2006.

[27] M. Gaetani-Liseo, C. Alonso, L. Seguier, and B. Jammes, 'Impact on Energy Saving of Active Phase Count Control to a DC/DC Converter in a DC Micro Grid', in *2018 7th International Conference on Renewable Energy Research and Applications (ICRERA)*, Paris, pp. 511–516, Oct. 2018.

[28] 'Datasheet YUASA', *NP VALVE REGULATED LEAD ACID BATTERY.* NP12-6, 6V 12Ah.

[29] W. E. Hart *et al.*, *Pyomo — Optimization Modeling in Python*, vol. 67. Cham: Springer International Publishing, 2017.

[30] G.-L. Margot, A. Corinne, and J. Bruno, 'Identification of ESS Degradations Related to Their Uses in Micro-Grids: Application to a Building Lighting Network with VRLA Batteries', *EJEE*, vol. 23, no. 6, pp. 455–466, Dec. 2021.

[31] REDDY, T. B., & LINDEN, D. (2011). *Linden's handbook of batteries*. New York, McGraw-Hill.

[32] D. U. Sauer, G. Bopp, A. Jossen, J. Garche, M. Rothert, and M. Wollny, 'State of charge - What do we really speak about ?', *International Communications Energy Conference*, Copenhague, Danemark, Jun. 1999.

[33] L. Zheng, J. Zhu, G. Wang, T. He, and Y. Wei, 'Novel methods for estimating lithium-ion battery state of energy and maximum available energy', *Applied Energy*, vol. 178, pp. 1–8, Sep. 2016.